\documentclass[twoside,11pt]{amsart}
\usepackage{amsmath,amssymb,amsthm}

\begin{document}

\setlength{\textwidth}{145mm} \setlength{\textheight}{203mm}
\setlength{\parindent}{0mm} \setlength{\parskip}{2pt plus 2pt}

\frenchspacing



\numberwithin{equation}{section}
\newtheorem{thm}{Theorem}[section]
\newtheorem{lem}[thm]{Lemma}
\newtheorem{prop}[thm]{Proposition}
\newtheorem{cor}[thm]{Corollary}
\newtheorem{probl}[thm]{Problem}

\newtheorem{defn}{Definition}[section]
\newtheorem{rem}{Remark}[section]
\newtheorem{exa}{Example}



\newcommand{\X}{\mathfrak{X}}
\newcommand{\B}{\mathcal{B}}
\newcommand{\s}{\mathfrak{S}}
\newcommand{\g}{\mathfrak{g}}
\newcommand{\W}{\mathcal{W}}
\newcommand{\Lgr}{\mathrm{L}}
\newcommand{\dd}{\mathrm{d}}

\newcommand{\pd}{\partial}
\newcommand{\ddx}{\frac{\pd}{\pd x^i}}
\newcommand{\ddy}{\frac{\pd}{\pd y^i}}
\newcommand{\ddu}{\frac{\pd}{\pd u^i}}
\newcommand{\ddv}{\frac{\pd}{\pd v^i}}

\newcommand{\diag}{\mathrm{diag}}
\newcommand{\End}{\mathrm{End}}
\newcommand{\im}{\mathrm{Im}}
\newcommand{\id}{\mathrm{id}}

\newcommand{\ie}{i.e.}
\newfont{\w}{msbm9 scaled\magstep1}
\def\R{\mbox{\w R}}
\newcommand{\norm}[1]{\left\Vert#1\right\Vert ^2}
\newcommand{\nnorm}[1]{\left\Vert#1\right\Vert ^{*2}}
\newcommand{\nN}{\norm{N}}
\newcommand{\nJ}{\norm{\nabla J}}
\newcommand{\nnJ}{\nnorm{\nabla J}}
\newcommand{\tr}{{\rm tr}}

\newcommand{\thmref}[1]{Theorem~\ref{#1}}
\newcommand{\propref}[1]{Proposition~\ref{#1}}
\newcommand{\secref}[1]{\S\ref{#1}}
\newcommand{\lemref}[1]{Lemma~\ref{#1}}
\newcommand{\dfnref}[1]{Definition~\ref{#1}}

\frenchspacing


\title[A connection with skew symmetric torsion ...]{A
connection with skew symmetric torsion and K\"ahler curvature
tensor on quasi-K\"ahler manifolds with Norden metric}

\author{Dimitar Mekerov}

\maketitle

{\small
\textbf{Abstract} \\
There is considered a connection with skew symmetric torsion on a
quasi-K\"ahler manifold with Norden metric. Some necessary and
sufficient conditions are derived for the corresponding curvature
tensor to be K\"ahlerian. In the case when this tensor is
K\"ahlerian, some relations are obtained between its scalar
curvature and the scalar curvature of other curvature tensors.
Conditions are given for the
considered manifolds to be isotropic-K\"ahler.}


\section{Preliminaries}

Let $(M,J,g)$ be a $2n$-dimensional \emph{almost complex manifold
with Norden metric}, i.e. $M$ is a differentiable manifold with an
almost complex structure $J$ and a metric $g$ such that
\begin{equation}\label{1.1}
J^2x=-x, \qquad g(Jx,Jy)=-g(x,y)
\end{equation}
for arbitrary $x$, $y$ of the algebra $\X(M)$ on the smooth vector
fields on $M$.

The associated metric $\tilde{g}$ of $g$ on $M$ is defined by
$\tilde{g}(x,y)=g(x,Jy)$. Both metrics are necessarily of
signature $(n,n)$. The manifold $(M,J,\tilde{g})$ is an almost
complex manifold with Norden metric, too.

Further, $x$, $y$, $z$, $w$ will stand for arbitrary elements of
$\X(M)$.

A classification of the almost complex manifolds with Norden
metric is given in \cite{GaBo}. This classification is made with
respect to the tensor field $F$ of type (0,3) defined by
\begin{equation}\label{1.2}
F(x,y,z)=g\bigl( \left( \nabla_x J \right)y,z\bigr),
\end{equation}
where $\nabla$ is the Levi-Civita connection of $g$. The tensor
$F$ has the following properties
\begin{equation}\label{1.3}
F(x,y,z)=F(x,z,y)=F(x,Jy,Jz).
\end{equation}

Among the basic classes $\W_1$, $\W_2$, $\W_3$ of this
classification, the almost complex structure is nonintegrable only
in the class $\W_3$. This is the class of the so-called
\emph{quasi-K\"ahler manifolds with Norden metric}, which we call
briefly \emph{$\W_3$-manifolds}. This class is characterized by
the condition
\begin{equation}\label{1.4}
\mathop{\s} \limits_{x,y,z} F(x,y,z)=0,
\end{equation}
where $\s $ is the cyclic sum by three arguments. The special
class $\W_0$ of the \emph{K\"ahler manifolds with Norden metric}
belonging to any other class is determined by the condition
$F(x,y,z)=0$.

Let $R$ be the curvature tensor of $\nabla$, i.e.
$R(x,y)z=\nabla_x \bigl(\nabla_y z\bigr) - \nabla_y \bigl(\nabla_x
z\bigr) - \nabla_{[x,y]}z$. The corresponding tensor of type
$(0,4)$ is determined by $R(x,y,z,w)$ $=g(R(x,y)z,w)$.

The following Ricci identity for almost complex manifolds with
Norden metric is known
\begin{equation}\label{1.5}
\bigl(\nabla_x F\bigr)(y,z,w)-\bigl(\nabla_y
F\bigr)(x,z,w)=R(x,y,Jz,w) - R(x,y,z,Jw).
\end{equation}

The components of the inverse matrix of $g$ are denoted by
$g^{ij}$ with respect to the basis $\{e_i\}$ of the tangent space
$T_pM$ of $M$ at a point $p\in M$.

The \emph{square norm of $\nabla J$} is defined by
\begin{equation}\label{1.6}
    \norm{\nabla J}=g^{ij}g^{ks}
    g\bigl(\left(\nabla_{e_i} J\right)e_k,\left(\nabla_{e_j}
    J\right)e_s\bigr).
\end{equation}
In \cite{MekMan-1} the following equation is proved for a
$\W_3$-manifold
\begin{equation}\label{1.7}
    \norm{\nabla J}=-2g^{ij}g^{ks}
    g\bigl(\left(\nabla_{e_i} J\right)e_k,\left(\nabla_{e_s}
    J\right)e_j\bigr).
\end{equation}

An almost complex manifold with Norden metric $(M,J,g)$ is
K\"ahlerian iff $\nabla J=0$. It is clear that we have
$\norm{\nabla J}=0$ for such a manifold, but the inverse one is
not always true. An almost complex manifold with Norden metric
with $\norm{\nabla J}=0$ is called an \emph{isotropic-K\"ahlerian}
in \cite{MekMan-1}.

The Ricci tensor $\rho$ for the curvature tensor $R$ and the
scalar curvature $\tau$ for $R$ are defined respectively by
\begin{equation}\label{1.8}
    \rho(x,y)=g^{ij}R(e_i,x,y,e_j),\qquad
    \tau=g^{ij}\rho(e_i,e_j),
\end{equation}
and their associated quantities $\rho^*$ and $\tau^*$ are
determined respectively by
\begin{equation}\label{1.9}
    \rho^*(x,y)=g^{ij}R(e_i,x,y,Je_j),\qquad
    \tau^*=g^{ij}\rho(e_i,Je_j).
\end{equation}
Similarly, the Ricci tensor and the scalar curvature are
determined for each \emph{curvature-like tensor (curvature
tensor)} $L$, i.e.
 for the tensor $L$ with the following properties:
\begin{equation}\label{1.10}
    L(x,y,z,w)=-L(y,x,z,w)=-L(x,y,w,z),
\end{equation}
\begin{equation}\label{1.11}
    \mathop{\s} \limits_{x,y,z} L(x,y,z,w)=0 \qquad \text{(first Bianchi
    identity)}.
\end{equation}
A curvature-like tensor is called a \emph{K\"ahler tensor} if it
has the property
\begin{equation}\label{1.12}
    L(x,y,Jz,Jw)=-L(x,y,z,w).
\end{equation}

The characteristic condition \eqref{1.4} for $\W_3$ is equivalent
to each of the following conditions \cite{MekMan-1}:
\begin{equation}\label{1.13}
    \mathop{\s} \limits_{x,y,z} F(Jx,y,z)=0,
\end{equation}
\begin{equation}\label{1.14}
    \left(\nabla_{x} J\right)Jy+\left(\nabla_{y} J\right)Jx
    +\left(\nabla_{Jx} J\right)y+\left(\nabla_{Jy} J\right)x=0.
\end{equation}

The following identity for a $\W_3$-manifold is known from
\cite{Me}:
\begin{equation}\label{1.15}
\begin{array}{l}
  \mathop{\s} \limits_{x,y,z} \bigl\{
  R(x,Jy,Jz,w)-R(x,Jy,z,Jw)\bigr.
  \\[4pt]
   \bigl.
  \phantom{\mathop{\s} \limits_{x,y,z} }
  +R(Jx,y,z,Jw)-R(Jx,y,Jz,w)\bigr\}
  \\[4pt]
  =-\mathop{\s}\limits_{x,y,z}
    g\Bigl(\bigl(\nabla_x J\bigr)y+\bigl(\nabla_y J\bigr)x,
     \bigl(\nabla_z J\bigr)w+\bigl(\nabla_w J\bigr)z\Bigr). \\[4pt]
\end{array}
\end{equation}

\section{A connection with skew symmetric torsion on a
$\W_3$-manifold}

A linear connection $\nabla'$ on an almost complex manifold with
Norden metric $(M,J,g)$ preserving $J$ and $g$, i.e.
$\nabla'J=\nabla'g=0$, is called a \emph{natural connection}
\cite{GaMi}. If $T$ is a torsion tensor of $\nabla'$, i.e.
$T(x,y)=\nabla'_x y-\nabla'_y x-[x, y]$, then the corresponding
tensor field of type (0,3) is determined by
$T(x,y,z)=g(T(x,y),z)$.

The connections with skew symmetric torsion are of particular
interest in the string theory \cite{Stro}. In mathematics this
connection was used by Bismut \cite{Bis} to prove the local index
theorem for non-K\"ahler Hermitian manifolds.

In this paper we consider a natural connection $\nabla'$  with
skew symmetric torsion on quasi-K\"ahler manifolds with Norden
metric whose curvature tensor has the properties of the curvature
tensor of a K\"ahler manifold with Norden metric. This connection
is determined by
\begin{equation}\label{1.16}
    \nabla'_{x} y=\nabla_{x} y+Q(x,y),
\end{equation}
where
\begin{equation}\label{1.17}
    Q(x,y)=\frac{1}{4}\Bigl\{\bigl(\nabla_x J\bigr)Jy-\bigl(\nabla_Jx
    J\bigr)y-2\bigl(\nabla_y J\bigr)Jx\Bigr\}.
\end{equation}
For the torsion tensor $T$ of $\nabla'$ we have $T(x,y)=2Q(x,y)$.
We denote
\begin{equation}\label{1.18}
    Q(y,z,w)=g(Q(y,z),w)
\end{equation}
and according to \eqref{1.2}, \eqref{1.3}, \eqref{1.17} and
\eqref{1.18} we obtain
\begin{equation}\label{1.19}
    Q(y,z,w)=-\frac{1}{4}\mathop{\s} \limits_{y,z,w} F(y,z,Jw).
\end{equation}


\section{Conditions for the curvature tensor of the connection $\nabla'$
on $\W_3$-manifolds to be K\"ahlerian}

Let $R'$ be the curvature tensor of the connection $\nabla'$ on a
$\W_3$-manifold $(M,J,g)$ determined by \eqref{1.16}, i.e.
\begin{equation}\label{2.1}
    R'(x,y)z=\nabla'_x \bigl(\nabla'_y z\bigr) - \nabla'_y \bigl(\nabla'_x
z\bigr) - \nabla'_{[x,y]}z.
\end{equation}
The corresponding tensor of type $(0,4)$ is determined by
$R'(x,y,z,w)=g(R'(x,y)z,w)$. According to \eqref{1.16},
\eqref{1.17} and \eqref{1.18}, we have
\begin{equation}\label{2.2}
    g\left(\nabla'_{x} y, z\right)=g\left(\nabla_{x} y,
    z\right)+Q(x,y,z).
\end{equation}
Since $\nabla g=\nabla' g=0$ then \eqref{2.1}, \eqref{2.2} and
\eqref{1.16} imply
\begin{equation}\label{2.3}
\begin{array}{l}
    R'(x,y,z,w)=R(x,y,z,w)+\left(\nabla_x Q\right)(y,z,w)-\left(\nabla_y
    Q\right)(x,z,w)
  \\[4pt]
\phantom{R'(x,y,z,w)=}
       -g\left(Q(y,z),Q(x,w)\right)+g\left(Q(x,z),Q(y,w)\right).
\end{array}
\end{equation}

The last equality implies the property \eqref{1.10} for $R'$ and
since $\nabla' J=0$ then \eqref{1.12} is valid, too. Therefore
$R'$ becomes K\"ahlerian if the condition \eqref{1.11} is
fulfilled for this tensor. Because of \eqref{2.3} the equality
\eqref{1.11} is valid for $R'$ iff
\begin{equation}\label{2.4}
\begin{array}{l}
    \mathop{\s} \limits_{x,y,z} \bigl\{
    \left(\nabla_x Q\right)(y,z,w)-\left(\nabla_y
    Q\right)(x,z,w)\bigr.
  \\[4pt]
\phantom{\mathop{\s} \limits_{x,y,z} \bigl\{\bigr.}\bigl.
       -g\left(Q(y,z),Q(x,w)\right)+g\left(Q(x,z),Q(y,w)\right)\bigr\}=0.
\end{array}
\end{equation}
Since $Q$ is a totally skew symmetric tensor then \eqref{2.4} gets
the form
\begin{equation}\label{2.5}
    \mathop{\s} \limits_{x,y,z} \bigl\{
    \bigl(\nabla_x Q\bigr)(y,z,w)\bigr\}
       =\mathop{\s} \limits_{x,y,z} \bigl\{g\bigl(Q(y,z),Q(x,w)\bigr)\bigr\}.
\end{equation}
The last equality implies immediately
\[
\begin{array}{l}
    \left(\nabla_x Q\right)(y,z,w)-\left(\nabla_y
    Q\right)(x,z,w)
  \\[4pt]
   = -\left(\nabla_z Q\right)(x,y,w)+\mathop{\s} \limits_{x,y,z} \bigl\{
   g\left(Q(x,y),Q(z,w)\right)\bigr\}
\end{array}
\]
and then \eqref{2.3} gets the form
\begin{equation}\label{2.6}
\begin{array}{l}
    R'(x,y,z,w)=R(x,y,z,w)
  \\[4pt]
\phantom{R'(x,y,z,w)=}
       -\left(\nabla_z Q\right)(x,y,w)+g\left(Q(x,y),Q(z,w)\right).
\end{array}
\end{equation}
In \eqref{2.6} we substitute $y\leftrightarrow w$ and we add the
obtained equality to \eqref{2.6}. Then we receive
\begin{equation}\label{2.7}
\begin{array}{l}
    R'(x,y,z,w)+R'(z,y,x,w)=R(x,y,z,w)+R(z,y,x,w)
  \\[4pt]
\phantom{R'(x,y,z,w)}
       +g\left(Q(x,y),Q(z,w)\right)+g\left(Q(z,y),Q(x,w)\right).
\end{array}
\end{equation}
Now we substitute $z\leftrightarrow w$ in \eqref{2.7} and  then we
subtract the obtained equality from \eqref{2.7}. Using the
properties of $R$ and $R'$ in the last equality we finally obtain
the following identity, which is equivalent to \eqref{2.4}:
\begin{equation}\label{2.8}
\begin{array}{l}
    3R'(x,y,z,w)=3R(x,y,z,w)+2g\left(Q(x,y),Q(z,w)\right)
  \\[4pt]
\phantom{3R'(x,y,z,w)=}
       +g\left(Q(z,y),Q(x,w)\right)+g\left(Q(x,z),Q(y,w)\right).
\end{array}
\end{equation}

In this way we proved the following
\begin{thm}\label{thm2.1}
Let $(M,J,g)$ be a $\W_3$-manifold and $\nabla'$ be the connection
determined by \eqref{1.16}. Then the curvature tensor $R'$ for
$\nabla'$ is K\"ahlerian iff the condition \eqref{2.8} is valid.
\end{thm}

Obviously the tensor $P$ defined by
\begin{equation}\label{2.9}
\begin{array}{l}
    P(x,y,z,w)=2g\left(Q(x,y),Q(z,w)\right)
  \\[4pt]
\phantom{P(x,y,z,w)=}
       +g\left(Q(z,y),Q(x,w)\right)+g\left(Q(x,z),Q(y,w)\right)
\end{array}
\end{equation}
satisfies the properties \eqref{1.10} and \eqref{1.11}, i.e. $P$
is a curvature-like tensor. Then from \thmref{thm2.1} we obtain
the following
\begin{cor}\label{cor2.2}
Let $(M,J,g)$ be a $\W_3$-manifold with a K\"ahler curvature
tensor $R'$ for the connection $\nabla'$ determined by
\eqref{1.16}. Then the tensor $P$ determined by \eqref{2.9} is
K\"ahlerian iff the curvature tensor $R$ is K\"ahlerian.
\end{cor}

Using \eqref{1.3}, \eqref{1.5}, \eqref{1.18}, \eqref{1.19},
\eqref{2.1} and the first Bianchi identity for $R$, we get the
following identity, which is equivalent to \eqref{2.5}:
\begin{equation}\label{2.10}
    \mathop{\s} \limits_{x,y,z} \bigl\{ \left(\nabla_w F\right)(x,z,Jy) \bigr\}=A(x,y,z,w),
\end{equation}
where
\begin{equation}\label{2.11}
\begin{array}{l}
    A(x,y,z,w)=\mathop{\s} \limits_{x,y,z} \bigl\{
    R(x,y,Jz,Jw)+R(Jx,Jy,z,w)
  \\[4pt]
\phantom{A(x,y,z,w)=\mathop{\s} \limits_{x,y,z}}
    +4g\left(Q(x,y),Q(z,w)\right)
       -g\bigl(\left(\nabla_x J\right)y, \left(\nabla_w
       J\right)z\bigr)
  \\[4pt]
\phantom{A(x,y,z,w)=\mathop{\s} \limits_{x,y,z}}
       +g\bigl(\left(\nabla_x J\right)y - \left(\nabla_y J\right)x, \left(\nabla_z J\right)w\bigr)
    \bigr\}.
\end{array}
\end{equation}
According to the properties of $F$, from \eqref{2.10} and
\eqref{2.11} we obtain
\begin{equation}\label{2.12}
A(Jx,y,z,w)+A(x,Jy,z,w)+A(x,y,Jz,w)-A(Jx,Jy,Jz,w)=0.
\end{equation}
Because of \eqref{2.11} the last equality implies
\begin{equation}\label{2.13}
\begin{array}{l}
   \mathop{\s} \limits_{x,y,z} \bigl\{
      g\bigl(\left(\nabla_x J\right)Jy+\left(\nabla_{Jx} J\right)y, \left(\nabla_w
       J\right)z+\left(\nabla_{Jz}J\right)Jw-\left(\nabla_z J\right)w\bigr)
    \bigr\}
  \\[4pt]
    =2\mathop{\s} \limits_{x,y,z} \bigl\{
    g\bigl(Q(x,y),Q(Jz,w)\bigr)+g\bigl(Q(Jx,y),Q(z,w)\bigr)
  \\[4pt]
\phantom{\mathop{\s} \limits_{x,y,z}=}
    +g\bigl(Q(x,Jy),Q(z,w)\bigr)-g\bigl(Q(Jx,Jy),Q(Jz,w)\bigr)\bigr\}.
\end{array}
\end{equation}
Having in mind $Q(x,Jy)=JQ(x,y)-\left(\nabla_x J\right)y$ and
\eqref{1.14}, from \eqref{2.13} we get the following identity,
equivalent to \eqref{2.5}
\begin{equation}\label{2.14}
   \mathop{\s} \limits_{x,y,z} \bigl\{
      g\bigl(\left(\nabla_x J\right)Jy+\left(\nabla_{Jx} J\right)y,
      \left(\nabla_z
       J\right)Jw+\left(\nabla_{Jz}J\right)w\bigr)
    \bigr\}=0.
\end{equation}
Then the following theorem is satisfied.
\begin{thm}\label{thm2.3}
Let $(M,J,g)$ be a $\W_3$-manifold and $\nabla'$ be the
connection determined by \eqref{1.16}. Then the curvature tensor
$R'$ for $\nabla'$ is K\"ahlerian iff the condition \eqref{2.14}
is valid.
\end{thm}

It is easy to verify that the properties \eqref{1.10},
\eqref{1.11} and \eqref{1.12} are valid for the tensor $H$ defined
by
\begin{equation}\label{2.15}
      H(x,y,z,w)=g\bigl(\left(\nabla_x J\right)Jy+\left(\nabla_{Jx} J\right)y,
      \left(\nabla_z
       J\right)Jw+\left(\nabla_{Jz}J\right)w\bigr).
\end{equation}
Then \thmref{thm2.3} implies the following
\begin{cor}\label{cor2.4}
Let $(M,J,g)$ be a $\W_3$-manifold and $\nabla'$ be the
connection determined by \eqref{1.16}. Then the curvature tensor
$R'$ for $\nabla'$ is K\"ah\-lerian iff the tensor $H$ determined
by \eqref{2.15} is K\"ahlerian.
\end{cor}


\section{Scalar curvatures on a $\W_3$-manifold with K\"ahler curvature tensor of the connection $\nabla'$}

Let $(M,J,g)$ be a $\W_3$-manifold with K\"ahler curvature tensor
of the connection and $\nabla'$ be determined by \eqref{1.16}.
Then the tensor $H$ determined by \eqref{2.15} is also K\"ahlerian
whereas the curvature tensor $R$ and the tensor $P$ determined by
\eqref{2.9} are curvature-like. We denote the scalar curvatures of
$R$, $R'$, $P$ and $H$ by $\tau$, $\tau'$, $\tau(P)$ and
$\tau(H)$, respectively, and their associated scalar curvatures by
$\tau^*$, $\tau'^*$, $\tau^*(P)$ and $\tau^*(H)$, respectively. We
denote  the associated square norm of $\nabla J$ with respect to
$\tilde{g}$ by $\nnJ$.

The equalities \eqref{2.8} and \eqref{2.9} imply immediately
\begin{equation}\label{3.1}
      3\tau'=3\tau+\tau(P),
\end{equation}
\begin{equation}\label{3.2}
      3\tau'^*=3\tau^*+\tau^*(P),
\end{equation}
\begin{equation}\label{3.3}
      \tau(P)=3g^{ij}g^{ks}g\bigl(Q(e_i,e_k),Q(e_s,e_j)\bigr).
\end{equation}
We obtain $g^{ij}F(e_i,e_j,z)=g^{ij}F(e_i,Je_j,z)=0$ from
\eqref{1.4}. The last equality and \eqref{1.19} imply
$g^{ij}Q(e_i,e_j)=0$. Then, having in mind \eqref{3.3}, we get
$\tau(P)=\frac{3}{8}\left(3\nJ+2\nnJ\right)$. Because of the
antisym\-metry of $Q$, \eqref{3.3} implies
$\tau(P)=\frac{3}{8}\left(3\nJ+\nnJ\right)$. In this way we obtain
\begin{equation}\label{3.4}
      \tau(P)=\frac{9}{8}\nJ.
\end{equation}
From \eqref{3.1} and  \eqref{3.4} we have
\begin{equation}\label{3.5}
      \tau'=\tau+\frac{3}{8}\nJ.
\end{equation}
By virtue of \eqref{2.9} we get
$\tau^*(P)=3g^{ij}g^{ks}g\bigl(Q(e_i,e_k),Q(Je_s,e_j)\bigr)$, from
where
\begin{equation}\label{3.6}
      \tau^*(P)=-\frac{3}{8}\nJ.
\end{equation}
Then, according to \eqref{1.15} and \eqref{3.2} we have
\begin{equation}\label{3.7}
      \tau'^*=\tau^*-\frac{1}{8}\nJ.
\end{equation}
The equalities \eqref{3.5} and \eqref{3.7}  imply
\begin{equation}\label{3.8}
      \tau'+3\tau'^*=\tau+3\tau^*.
\end{equation}
Using \eqref{1.14} and \eqref{2.15} we obtain
\begin{equation}\label{3.9}
      \tau(H)=\tau^*(H)=2\nJ.
\end{equation}
Then, from  \eqref{3.5},  \eqref{3.7} and \eqref{3.9} the
following equalities are valid
\begin{equation}\label{3.10}
      \tau'=\tau+\frac{3}{16}\tau(H),
\end{equation}
\begin{equation}\label{3.11}
      \tau'^*=\tau^*-\frac{1}{16}\tau(H).
\end{equation}
By virtue of \eqref{3.4}, \eqref{3.5}, \eqref{3.6}, \eqref{3.7}
and \eqref{3.9}, we get the following
\begin{thm}\label{thm-3.1}
Let $(M,J,g)$ be a $\W_3$-manifold with K\"ahler curvature tensor
$R'$ of the connection $\nabla'$ determined by \eqref{1.16}. Then
$(M,J,g)$ is an isotropic-K\"ahler  manifold iff an arbitrary one
of the quantities $\tau-\tau'$, $\tau^*-\tau'^*$, $\tau(P)$,
$\tau^*(P)$, $\tau(H)$, $\tau^*(H)$ is zero.
\end{thm}

Now, let $(M,J,g)$ be a 4-dimensional $\W_3$-manifold. Since $R'$
is a K\"ahler tensor, according to \cite{Teof-2} we have
\begin{equation}\label{3.12}
      R'=\nu'(\pi_1-\pi_2)+\nu'^*\pi_3,
\end{equation}
where $\nu'=\frac{\tau'}{8}$, $\nu'^*=\frac{\tau'^*}{8}$ and
\[
\begin{array}{l}
\pi_1(x,y,z,w)=g(y,z)g(x,w)-g(x,z)g(y,w),\\[4pt]
\pi_2(x,y,z,w)=g(y,Jz)g(x,Jw)-g(x,Jz)g(y,Jw),\\[4pt]
\pi_3(x,y,z,w)=-g(y,z)g(x,Jw)+g(x,z)g(y,Jw),\\[4pt]
\phantom{\pi_2(x,y,z,w)=} -g(y,Jz)g(x,w)+g(x,Jz)g(y,w).
\end{array}
\]
According to \eqref{3.5}, \eqref{3.7}, \eqref{3.12} and
\eqref{2.9}, from \eqref{2.8} we obtain
\begin{equation}\label{3.13}
      R=\frac{1}{8}\left\{\left(\tau+\frac{3}{8}\nJ\right)\left(\pi_1-\pi_2\right)
      +\left(\tau^*-\frac{1}{8}\nJ\right)\pi_3\right\}-\frac{1}{3}P.
\end{equation}
Then we have the following
\begin{thm}\label{thm-3.2}
Let $(M,J,g)$ be a 4-dimensional $\W_3$-manifold with K\"ahler
curvature tensor $R'$ of the connection $\nabla'$ determined by
\eqref{1.16}. Then $(M,J,g)$ is an isotropic-K\"ahler  manifold
iff
\[
R=\frac{1}{8}\left\{\tau\left(\pi_1-\pi_2\right)
      +\tau^*\pi_3\right\}-\frac{1}{3}P.
\]
\end{thm}

Because of \eqref{3.9} and \eqref{3.13} the following theorem is
valid.
\begin{thm}\label{thm-3.3}
Let $(M,J,g)$ be a 4-dimensional $\W_3$-manifold with K\"ahler
curvature tensor $R'$ of the connection $\nabla'$ determined by
\eqref{1.16}. Then we have
\[
      R=\frac{1}{128}\bigl\{\left(16\tau+\tau(H)\right)\left(\pi_1-\pi_2\right)
      +\left(16\tau^*-\tau(H)\right)\pi_3\bigr\}-\frac{1}{3}P.
\]
\end{thm}


\bigskip

\textit{Dimitar Mekerov\\
University of Plovdiv\\
Faculty of Mathematics and Informatics
\\
Department of Geometry\\
236 Bulgaria Blvd.\\
Plovdiv 4003\\
Bulgaria
\\
e-mail: mircho@uni-plovdiv.bg}

\end{document}